\documentclass[11pt,a4paper]{article}
\usepackage{a4wide,amsfonts,amsmath,latexsym,amssymb,euscript,eufrak,graphicx,units,mathrsfs}

\usepackage{stmaryrd}
\usepackage{float}
\newfloat{figure}{H}{lof}
\floatname{figure}{\figurename}

\DeclareMathAlphabet{\eufrak}{U}{}{}{}  
\SetMathAlphabet\eufrak{normal}{U}{euf}{m}{n}
\SetMathAlphabet\eufrak{bold}{U}{euf}{b}{n}

\usepackage{amsmath, amsthm, amsfonts, amssymb}

\numberwithin{equation}{section}

\newcommand{\diff}[1]{\operatorname{d}\ifthenelse{\equal{#1}{}}{\,}{#1}}
\def\real{{\mathord{{\rm I\kern-2.8pt R}}}}        
\def\inte{{\mathord{{\rm I\kern-2.8pt N}}}}
\def\PP{{\mathord{{\rm I\kern-2.8pt P}}}}

\def\real{{\mathord{\mathbb R}}}

\def\inte{{\mathord{\mathbb N}}}
\def\Z{{\mathord{\mathbb Z}}}
\def\D{{\mathord{\mathbb D}}}

\def\R{\right}

\def\R{\mathbb{R}}

\def\E{\mathbb{E}}

\def\real{\mathbb{R}}

\newtheorem{prop}{Proposition}[section]

\newtheorem{lemma}[prop]{Lemma}

\newtheorem{theorem}[prop]{Theorem}

\allowdisplaybreaks

\begin{document}

\begin{center}
{\large \textbf{The functional Breuer-Major theorem}}\\[0pt]
~\\[0pt]
 Ivan Nourdin
\footnote{%
Universit\'e du Luxembourg, Maison du Nombre, 6 avenue de la Fonte, 
L-4364 Esch-sur-Alzette, Grand Duchy of Luxembourg.
{\tt ivan.nourdin@uni.lu}}
and David Nualart  
\footnote{%
Department of Mathematics, University of Kansas, Lawrence, KS 66045, USA.
{\tt nualart@ku.edu}  }
 \footnote{David Nualart was supported by the NSF grant  DMS 1811181}
\\[0pt]
{\it Universit\'e du Luxembourg and University of Kansas}\\
~\\[0pt]
\end{center}

{\small \noindent \textbf{Abstract:} 
Let $X=\{ X_n\}_{n\in \mathbb{Z}}$  be  zero-mean stationary Gaussian  sequence of
random variables  with covariance function $\rho$ satisfying $\rho(0)=1$. Let  $\varphi:\mathbb{R}\to\mathbb{R}$ be a function such that $\E[\varphi(X_0)^2]<\infty$ 
and  assume that $\varphi$ has
Hermite rank $d \geq 1$.
The celebrated Breuer-Major theorem asserts that, if $\sum_{r\in\mathbb{Z}}
|\rho(r)|^d<\infty$ then
the finite dimensional distributions of $\frac1{\sqrt{n}}\sum_{i=0}^{\lfloor n\cdot\rfloor-1} \varphi(X_i)$ converge to those of $\sigma\,W$, where $W$ is a standard Brownian motion and $\sigma$ is some (explicit) constant.
Surprisingly, and despite the fact this theorem
has become over the years a prominent tool in a bunch of different areas,
a necessary and sufficient condition implying the weak convergence in the space ${\bf D}([0,1])$ of c\`adl\`ag functions
endowed with the Skorohod topology is still missing.
Our main goal in this paper is to fill this gap.
More precisely, by using suitable boundedness properties satisfied by the generator of the Ornstein-Uhlenbeck semigroup, we show that tightness holds under
the sufficient (and almost necessary) natural condition that $\E[|\varphi(X_0)|^{p}]<\infty$ for some $p>2$.
\normalsize
}

\section{Introduction}
Consider a zero-mean stationary Gaussian  sequence of random variables   $X=\{ X_n\}_{n\in \Z}$   
with covariance function 
$\E[X_nX_m ] =\rho(|n-m|)$ such that $\rho(0)=1$.  
Let $\gamma =N(0,1)$ be the standard Gaussian measure on $\R$. Consider a function
$\varphi \in L^2(\R,\gamma)$   of
Hermite rank $d \geq 1$, that is, $\varphi$ has a
series expansion given by
\begin{equation}\label{hermite}
  \varphi(x) = \sum_{q=d}^{\infty} c_q H_q(x), \,\,\, c_d\neq 0,
\end{equation}
where $H_q(x)$ is the $q$th Hermite polynomial with leading coefficient 1.

A classical central limit theorem, proved by Breuer and Major in \cite{BM}, asserts that under the condition
 \begin{equation} \label{h1}
\sum_{k \in \Z} |\rho(k)|^d < \infty,
\end{equation}
the \emph{finite-dimensional distributions}  of the process
  \begin{equation}\label{yn}
	  Y_n(t):= \frac{1}{\sqrt{n}} \sum_{i=0}^{ \lfloor nt \rfloor-1} \varphi(X_i) , \quad t\in[0,1]
	  \end{equation}
 converge to those of $\sigma W   $  as $n$ tends to infinity, where $W=\{W_t\}_{t\in [0,1]}$ is a standard Brownian motion and
 \begin{equation}\label{bm.sig}
		\sigma^2 = \sum_{q=d}^\infty q! c_q^2 \sum_{k \in \Z} \rho(k)^q.
	\end{equation}
	Observe that $|\rho(k)|=|\E[X_k X_0]|\leq \rho(0)=1$ by Cauchy-Schwarz, and thus $\sigma^2$ is well defined under the integrabilility assumption (\ref{h1}) imposed on $\rho$.
We also refer the reader to \cite[Chapter 7]{IvanGioBook}, where a modern proof of the Breuer-Major theorem is given, by means of the recent Malliavin-Stein approach.

What about the \emph{functional convergence}, that is, convergence in law of $Y_n$ to $\sigma W$ in  the  space ${\bf D}([0,1])$ endowed with the Skorohod topology?
The best-to-date available criterion ensuring tightness for $Y_n$ is due to Ben Hariz \cite{BenHariz}
and Chambers and Slud \cite{ChambersSlud} (the former being only a slight improvement with respect to the latter\footnote{Chambers and Slud criterion corresponds to Ben Hariz criterion (\ref{criterebenhariz}) with $R=\frac32$ and without the terms $\sum_{k\in\Z}|\rho(k)|^q$ all bounded by (\ref{h1}).}),  in the simpler situation where sums are replaced by integrals and convergences are understood in the space ${\bf C}([0,1])$ 
of continuous functions endowed with the uniform topology.
Transformed into our setting, the criterion in \cite{BenHariz,ChambersSlud} reads as follows\footnote{Compared to \cite{BenHariz}, condition (\ref{criterebenhariz}) is stated here with $\sqrt{q!}$ instead of  $(\sqrt{q!})^{-1}$   (since we work here with Hermite polynomials with leading coefficient 1) and with sums replacing integrals (since we work here in a discrete framework).}: tightness holds provided there exists $R>1$ such that
\begin{equation}\label{criterebenhariz}
\sum_{q=d}^\infty \sqrt{q!}|c_q|\left(
\sum_{k\in\Z}|\rho(k)|^q
\right)^\frac12R^q<\infty.
\end{equation}
But, in our opinion, condition (\ref{criterebenhariz}) is not meaningfull, for at least three reasons: $(i)$ it is not very natural, $(ii)$ it
is far from being optimal, and $(iii)$ it may be difficult to check it in practice, especially when the computation
of the Hermite coefficients $c_q$ appears to be tricky or even impossible.
Moreover, the proof given in \cite{BenHariz,ChambersSlud}   of the fact that  (\ref{criterebenhariz}) implies tightness can be simplified a lot, by proceeding as follows.  
Let us first recall that tightness in ${\bf D}([0,1])$
holds if their exist $p>2$ and $c>0$ such that, for all $n$,
\begin{equation}\label{billi}
\|Y_n(t)-Y_n(s)\|_{L^{p}(\Omega)}\leq c \left( \frac{ \lfloor nt \rfloor - \lfloor ns \rfloor } n \right) ^{1/2},\quad 0\leq s\leq t\leq 1
\end{equation}
(see Lemma \ref{billylemma} below).
Here, we have
\begin{eqnarray}
\|Y_n(t)-Y_n(s)\|_{L^{p}(\Omega)}&=&
\left\| \sum_{q=d}^{\infty} c_q \frac{1}{\sqrt{n}} \sum_{i=\lfloor ns\rfloor}^{ \lfloor nt \rfloor-1} H_q(X_i)\right\|_{L^{p}(\Omega)}\notag\\
&\leq&
 \sum_{q=d}^{\infty} |c_q| \left\|\frac{1}{\sqrt{n}} \sum_{i=\lfloor ns\rfloor}^{ \lfloor nt \rfloor-1} H_q(X_i)\right\|_{L^{p}(\Omega)}.\label{ineq1}
\end{eqnarray}
At this stage, a crucial observation is that $\sum_{i=\lfloor ns\rfloor}^{ \lfloor nt \rfloor-1} H_q(X_i)$ belongs to the $q$th Wiener chaos, where all $L^p(\Omega)$-norms are equivalent
by hypercontractivity.
More precisely,
\begin{equation}\label{hyper}
\left\|\frac{1}{\sqrt{n}} \sum_{i=\lfloor ns\rfloor}^{ \lfloor nt \rfloor-1} H_q(X_i)\right\|_{L^{p}(\Omega)}\leq (p-1)^{\frac{q}2}
\left\|\frac{1}{\sqrt{n}} \sum_{i=\lfloor ns\rfloor}^{ \lfloor nt \rfloor-1} H_q(X_i)\right\|_{L^2(\Omega)},
\end{equation}
see, e.g., \cite[Corollary 2.8.14]{IvanGioBook}.
The interest of the right-hand side of (\ref{hyper}) with respect to the left-hand side is that the former is straightforward to calculate and to estimate, as follows:
$$
\left\|\frac{1}{\sqrt{n}} \sum_{i=\lfloor ns\rfloor}^{ \lfloor nt \rfloor-1} H_q(X_i)\right\|_{L^2(\Omega)}^2
\leq 
\frac{ \lfloor nt \rfloor - \lfloor ns \rfloor } n \,\,
q!\sum_{k\in\Z}|\rho(k)|^q
.$$
By plugging this into (\ref{hyper}) and then into (\ref{ineq1}),
we obtain
$$
\|Y_n(t)-Y_n(s)\|_{L^p(\Omega)}\leq \left( \frac{ \lfloor nt \rfloor - \lfloor ns \rfloor } n \right) ^{1/2}\, \sum_{q=d}^{\infty} |c_q|
(p-1)^\frac{q}2\sqrt{q!}\left(\sum_{k\in\Z}|\rho(k)|^q\right)^\frac12
,
$$
implying in turn that (\ref{billi}) is satisfied (and then tightness) under (\ref{criterebenhariz}) with $R=\sqrt{p-1}>1$.

As we have just seen, the criterion (\ref{criterebenhariz}) of \cite{BenHariz,ChambersSlud} for tightness
is actually not so difficult  to prove. But on the other hand it is neither natural, nor easy to check in practice.
The main objective of this note is thus to provide a simpler sufficient condition for the convergence $Y_n\Rightarrow\sigma W$ to hold in  law in ${\bf D}([0,1])$ endowed with the  Skorohod topology.
Actually, our finding is that only a little more integrability of the function $\varphi$ is needed.

\begin{theorem} \label{thm1}
Let $X=\{ X_n \}_{n\in \mathbb{Z} } $ be a zero-mean Gaussian stationary  sequence with covariance function 
$\E[X_nX_m ] =\rho(|n-m|)$ such that $\rho(0)=1$. Consider a function  $\varphi \in L^2(\R,\gamma)$ with expansion (\ref{hermite}) and Hermite rank $d\ge 1$, and suppose that  $\sum_{k\in \mathbb{Z}}  |\rho(k) |^d < \infty $.
Finally, recall $Y_n$ from (\ref{yn}) , let $W=\{W_t\}_{t\in [0,1]}$ be a Brownian motion and let $\sigma ^2$ be defined in (\ref{bm.sig}). Then, as $n\to\infty$,
\begin{enumerate}
\item The finite-dimensional distributions of $Y_n$ converge to those of $\sigma W$;
\item If  $\varphi \in L^{p}(\R,\gamma)$  for some $p>2$, then
$Y_n$ converges in   law to $\sigma W$ in ${\bf D}([0,1])$ 
endowed with the Skorohod topology.
\end{enumerate}
\end{theorem}

We can prove a similar result in the space ${\bf C}([0,1])$ of continuous functions
endowed with the uniform topology. Of course, in this case we have to consider the linear interpolation $Z_n$ instead of $Y_n$, defined
as follows:
\begin{equation}\label{zn}
Z_n(t)=\frac{nt-\lfloor nt\rfloor}{\sqrt{n}}\,\varphi(X_{\lfloor nt\rfloor })+\frac{1}{\sqrt{n}} \sum_{i=0}^{ \lfloor nt \rfloor-1} \varphi(X_i),\quad t\in[0,1].
\end{equation}
\begin{theorem} \label{thm2}
Let $X=\{ X_n\}_{n\in \mathbb{Z} } $ be a zero-mean Gaussian stationary  sequence with covariance function 
$\E[X_nX_m ] =\rho(|n-m|)$ such that $\rho(0)=1$. Consider a function  $\varphi \in L^2(\R,\gamma)$ with expansion (\ref{hermite}) and Hermite rank $d\ge 1$, and suppose that  $\sum_{k\in \mathbb{Z}}  |\rho(k) |^d < \infty $.
Finally, recall $Z_n$ from (\ref{zn}), let $W=\{W_t\}_{t\in [0,1]}$ be a Brownian motion and let $\sigma ^2$ be defined in (\ref{bm.sig}). Then, as $n\to\infty$,
\begin{enumerate}
\item The finite-dimensional distributions of $Z_n$ converge to those of $\sigma W$;
\item If  $\varphi \in L^{p}(\R,\gamma)$  for some $p>2$, then
$Z_n$ converges in   law to $\sigma W$ in ${\bf C}([0,1])$ 
endowed with the uniform topology.
\end{enumerate}
\end{theorem}

The proof of Theorems \ref{thm1} and \ref{thm2} are
based on the  application of the techniques of Malliavin calculus,
especially Meyer inequalities from \cite{Meyer} (in the modern form taken from \cite{DavidBook}).
The method we employ is based on the representation 
$\varphi(X_i)=\delta^d\big((D(-L)^{-1})^d(\varphi(X_i))\big)$ where $\delta$, $D$ and $L$ are the usual Malliavin operators (see Section 2).
It
is robust enough to be used for other families of interest
than $Y_n$ and $Z_n$, see indeed \cite{JN} for an application to the self-intersection local time of the fractional Brownian motion, or Section 4 in the present paper for an extension of Theorem \ref{thm1} in a critical situation where $\sum_{|k|\leq n} |\rho(k)|^{d}$ diverges slowly enough when $n\to\infty$.

The rest of the paper is organized as follows. 
Section 2 contains some useful preliminaries on Malliavin calculus,
as well as some boundedness properties of the so-called shift operator, which is
our main tool in this paper.
The proof of Theorem \ref{thm1} (resp. \ref{thm2}) is given in Section 3 (resp. 4).
Finally, in Section 4 we provide an extension of Theorem \ref{thm1}
in the case where $\sum_{k\in\Z}|\rho(k)|^d$ explodes slowly.


\section{Preliminaries}  \label{sec2}
In this section, we gather several preliminary results that are needed for the proofs of  the main results of this paper. 

\subsection{Elements of Malliavin calculus with respect to the Wiener process}\label{sec1.1}

We refer the reader to the references \cite{IvanGioBook,DavidBook,DavidEulaliaBook} for a detailed account on the Malliavin calculus.  In this paper we will make use of the following notation and results.

First, let us introduce a specific realization of the sequence $\{X_k\}_{k\in\Z}$. The space
\[
\mathcal{H}:=\overline{{\rm span}\{X_k,\,k\in\Z\}}^{L^2(\Omega)}
\]
being a real separable Hilbert space, it is isometrically isomorphic to either $\R^N$ (for some $N\geq 1$) or
$L^2(\R_+)$. In both cases, there exists an isometry $\Phi:\mathcal{H}\to L^2(\R_+)$. Set $e_k=\Phi(X_k)$ for each $k\in\Z$. We have
\begin{equation}
\rho(k-l)=\E[X_kX_l]=\int_0^\infty e_k(x)e_l(x)dx,\quad k,l\geq 1.\label{ekrho}
\end{equation}
Let
$W=\{W(h),\,h\in L^2(\R_+)\}$ be the standard Wiener process, that is, a centered Gaussian family satisfying
$\E[W(h)W(g)]=\langle h,g\rangle_{L^2(\R_+)}$ for all $h,g\in L^2(\R_+)$.
We deduce immediately from (\ref{ekrho}) that
\[
\{X_k\}_{k\in\Z} \overset{\rm law}{=} \left\{W(e_k)\right\}_{k\in\Z}.
\]
Since, in this paper, the quantities we are interested in only depend on the law, starting from now and without loss of generality, we set
\begin{equation}\label{xk}
X_k:=W(e_k),\quad k\in\Z.
\end{equation}

 For integers $q\geq 1$,
 the $q$th Wiener chaos is the closed linear subspace
of $L^2(\Omega)$ that is generated by the random variables $\{H_q(W(h)),\,h\in L^2(\R_+),\,\|h\|_{L^2(\R_+)}=1\}$, where $H_q$ stands for the $q$th Hermite polynomial
defined by
\[
H_q(x)=(-1)^qe^{\frac{x^2}{2}}\frac{d^q}{dx^q}e^{-\frac{x^2}{2}},\quad q \geq 1,
\]
and $H_0(x)=1$.
For $q\geq 1$, it is known that the map 
\begin{equation}\label{rel}
I_q(h^{\otimes q})=H_q(W(h)),\quad h\in L^2(\R_+),\,\,\|h\|_{L^2(\R_+)}=1,
\end{equation} 
provides a linear isometry between  the set of symmetric square integrable functions $L^2_s(\R_+^q)$ (equipped with the modified 
norm $\sqrt{q!}\|\cdot\|_{L^2(\R_+^q)}$) and the $q$th Wiener chaos. By convention, $I_0(x)=x$ for all $x\in\R$.

It is well-known that any $F\in L^2(\Omega)$ measurable with respect to $W$ can be decomposed into Wiener chaos as follows:
\begin{equation}\label{decompo}
F=\E[F]+\sum_{q=1}^\infty I_q(f_q),
\end{equation}
where the kernels $f_q\in L^2_s(\R_+^q)$ are uniquely determined by $F$.

For a smooth and cylindrical random variable $F= f(W(h_1), \dots , W(h_n))$, with $h_i \in  L^2(\R_+)$ and $f \in C_b^{\infty}(\mathbb{R}^n)$ ($f$ and of its partial derivatives are bounded), we define its Malliavin derivative $D$ as the $L^2(\R_+)$-valued random variable given by
\[
DF = \sum_{i=1}^n \frac{\partial f}{\partial x_i} (W(h_1), \dots, W(h_n))h_i.
\]
By iteration, one can define the $k$-th derivative $D^k F$  as an element of $L^2(\Omega; L^2(\R_+^k))$. For any natural number $k$ and any real number $ p \geq 1$, we define  the Sobolev space $\mathbb{D}^{k,p}$  as the closure of the space of smooth and cylindrical random variables with respect to the norm $\|\cdot\|_{k,p}$ defined by 
\[
 \|F\|^p_{k,p} = \mathbb{E}(|F|^p) + \sum_{l=1}^k \mathbb{E}(\|D^l F\|^p_{L^2(\R_+^l)}).
\]
For any Hilbert space $V$ we denote by $\D^{k,p} (V)$ the corresponding space of $V$-valued random variables.

The divergence operator $\delta$ is defined as the adjoint of the derivative operator $D$.  An element $u \in L^2(\Omega; L^2(\R_+))$ belongs to the domain of $\delta$, denoted by ${\rm Dom}\, \delta$, if there is a constant $c_u$ depending on $u$ such that 
\[
|\mathbb{E} (\langle DF, u \rangle_{L^2(\R_+)})| \leq c_u \|F\|_{L^2(\Omega)}
\] for any $F \in \mathbb{D}^{1,2}$.  If $u \in {\rm Dom} \,\delta$, then the random variable $\delta(u)$ is defined by the duality relationship 
\begin{equation} \label{dua}
\mathbb{E}[F\delta(u)] = \mathbb{E}[\langle DF, u \rangle_{L^2(\R_+)}] \, ,
\end{equation}
which holds for any $F \in \mathbb{D}^{1,2}$.  
In a similar way we can introduce the iterated divergence operator $\delta^k$ for each integer $k\ge 2$, defined by the duality relationship 
\begin{equation} \label{dua2}
\mathbb{E}[F\delta^k(u)] = \mathbb{E}  \left[\langle D^kF, u \rangle_{L^2(\R_+^k)} \right] \, ,
\end{equation}
for any $F \in \mathbb{D}^{k,2}$, where $u\in  {\rm Dom}\, \delta^k \subset L^2(\Omega; L^2(\R_+^k))$.
If $u\in L^2_s(\R_+^k)$ is deterministic, then
\begin{equation} \label{equ44}
\delta^k(u)= I_k(u).
\end{equation}
For any $p>1$ and any integer $k\ge 1$, the operator $\delta^k$ is continuous  from $\D^{k,p} (L^2(\R_+^k))$ into  $L^p(\Omega)$, and we have the  inequality (see, for instance, \cite[Proposition 1.5.4]{DavidBook})
\begin{equation} \label{meyer1}
\|\delta ^k(v) \| _{L^p(\Omega)} \le c_p \sum_{j=0}^k  \| D^jv\|_{L^p(\Omega; L^2(\R_+^j))},
\end{equation}
for any $v\in \D^{k,p} (L^2(\R_+^k))$. This inequality is a consequence of Meyer inequalities (from \cite{Meyer}), which  states the equivalence in $L^p(\Omega)$, for any $p>1$, of the operators $D$ and $(-L) ^{1/2}$,
where $L$ is the infinitesimal generator of the Ornstein-Uhlenbeck semigroup $(P_t)_{t \geq 0}$ in $L^2(\Omega)$ defined as
$$P_tF=\sum_{q=0}^\infty e^{-qt} I_q(f_q),\,\,t\geq 0,\quad\mbox{and}\quad (-L)^r F=\sum_{q=1}^\infty q^r I_q(f_q),\,\,\, r\in\R,$$ if $F$ is given by  (\ref{decompo}). More precisely, their exist two constants $c_{i,p}$, $i=1,2$, such that, for any $F\in \D^{1,p}$,
\begin{equation} \label{meyer2}
c_{1,p}  \| DF\| _{L^p(\Omega, L^2(\R_+))} \le \| (-L)^{1/2} F\| _{ L^p(\Omega)} \le  c_{2,p}  \| DF\| _{L^p(\Omega, L^2(\R_+))}.
\end{equation}
More generally, we can state Meyer's inequalities in the general case (see \cite[Theorem 1.5.1]{DavidBook}): for any $p>1$ and any integer $k\geq 1$,  their exist two constants $c_{i,p,k}$, $i=1,2$, such that, for any $F\in \D^{1,p}$,
\begin{equation} \label{meyer3}
c_{1,k,p}  \| D^kF\| _{L^p(\Omega, L^2(\R^k_+))} \le \| (- L)^{k/2} F\| _{ L^p(\Omega)} \le  c_{2,k,p} \big( \| D^kF\| _{L^p(\Omega, L^2(\R^k_+))}
+\|F\|_{L^p(\Omega)}\big).
\end{equation}

\subsection{The shift operator}
Let  $\varphi\in L^2(\R,\gamma)$  be a function of Hermite rank $d\ge 1$ and expansion  (\ref{hermite}).
Consider  the function $\varphi_d$ defined by a shift of $d$ units in the coefficients, that is,
\begin{equation} \label{fd3}
\varphi_d= \sum_{q=d} ^\infty c_q H_{q-d}.
\end{equation}
It is immediately checked that $\varphi_d \in L^2(\R,\gamma)$.

Given (\ref{fd3})  and the relation (\ref{rel}) between Hermite polynomials and multiple stochastic integrals, the random variable $\varphi_d(W(h))$ 
admits the following
chaotic decomposition  when $h\in L^2(\R_+)$ has norm 1 :
$$
\varphi_d(W(h))=\sum_{q=d}^\infty c_q I_{q-d}(h^{\otimes(q-d)}).
$$
Moreover, we claim that $\varphi_d(W(h))$ belongs to 
$\mathbb{D}^{2,d}$.
Indeed, for any $k=1,\ldots,d$ we have
that
\[
D^k(\varphi_d(W(h))) = \sum_{q=d} ^\infty c_q (q-d) (q-d-1) \cdots (q-d-k+1) I_{q-d-k}(h^{\otimes(q-d-k)})h^{\otimes k},
\]
and this series converges in $L^2(\Omega,L^2(\R_+^k))$ since
\begin{align*}
\E\| D^k(\varphi_d(W(h))) \|^2_{L^2(\R_+^k)}
& = \sum_{q=d}^ \infty c_q^2  (q-d) ^2(q-d-1)^2\cdots (q-d-k+1)^2 (q-d-k)!  \\
&\leq  \sum_{q=d} ^\infty c_q^2  q!< \infty.
\end{align*}

\bigskip

The following two lemmas will play a crucial role in the sequel.

\begin{lemma}\label{crucial2}
Suppose that  $\varphi\in L^{2}(\R,\gamma)$ given by (\ref{hermite}) has Hermite rank $d\geq 1$ .
We have, for any $h\in L^2(\R_+)$ of norm 1,
\begin{eqnarray}\label{varphir}
\varphi(W(h))&=&\delta^d(\varphi_d(W(h))h^{\otimes d})\\
\varphi_d(W(h))h^{\otimes d}&=&(D (-L)^{-1})^d (\varphi(W(h)))\label{varphir2}\\
\label{varphid}
\varphi_d(W(h))&=&\langle (D (-L)^{-1})^d (\varphi(W(h))),h^{\otimes d}\rangle_{L^2(\R_+^d)}.
\end{eqnarray}
\end{lemma}
\noindent
{\it Proof}.
Using (\ref{equ44}) and  the relation (\ref{rel}) between Hermite polynomials and multiple stochastic integrals, we can write
\begin{align*}
\varphi(W(h))&=\sum_{q=d} ^\infty c_q H_{q}(W(h)) =\sum_{q=d} ^\infty c_q   I_q ( h^{\otimes q} )=\sum_{q=d} ^\infty c_q  \delta^q ( h^{\otimes q} )\\
&
=
\sum_{q=d} ^\infty c_q  \delta^d \left(  \delta^{q-d} \left(
 h^{\otimes q-d} \right)  h^{\otimes d}
\right) 
= \delta^d \left( \sum_{q=d} ^\infty c_q   I_{q-d} \left(
 h^{\otimes q-d} \right)     h^{\otimes d}
\right) \\
&= \delta^d \left( \sum_{q=d} ^\infty c_q   H_{q-d} (W(h))   h^{\otimes d}
\right)  =\delta^d(\varphi_d(W(h))h^{\otimes d}),
\end{align*}
which is (\ref{varphir}).
On the other hand, we can compute that
$$
(D(-L)^{-1})(\varphi(W(h)))=
\sum_{q=d}^\infty c_q\,I_{q-1}(h^{\otimes q-1})\,h.
$$
By iteration, we get
$$
(D(-L)^{-1})^d(\varphi(W(h)))=
\sum_{q=d}^\infty c_q\,I_{q-d}(h^{\otimes q-d})\,h^{\otimes d}
=
\sum_{q=d}^\infty c_q\,H_{q-d}(W(h))\,h^{\otimes d},
$$
and the desired conclusions (\ref{varphir2}) and then (\ref{varphid}) follow.\qed

\bigskip

\begin{lemma}\label{lmcrucial}
Suppose that  $\varphi\in L^{2}(\R,\gamma)$ given by (\ref{hermite}) has Hermite rank $d\geq 1$ and is such that 
$\E\big[|\varphi(N)|^{p}\big]<\infty$ for some $p>2$ and $N\sim N(0,1)$.
Then, for any $0\leq k\leq r\leq d$,
$$\sup_{\|h\|=1} \E\big[\| D^k(D (-L)^{-1})^r(\varphi(W(h))) \|^{p}_{L^2(\R_+^{r+k})}\big]<\infty,$$
where the supremum runs over the set of all square integrable functions $h\in L^2(\R_+)$ of norm 1.
\end{lemma}
\noindent
{\it Proof}.
The proof is by induction on $r$. When $r=0$, one has $k=0$ and $D^0(D (-L)^{-1})^0$ is the identity operator, so there is nothing to prove.

Suppose now that the conclusion of Lemma \ref{lmcrucial} holds true for some $r-1\in\{0,\ldots,d-1\}$, and
let us prove that it holds true for $r+1$ as well.

If $k=0$, we have $D^0(D (-L)^{-1})^r=(D (-L)^{-1})^r$.
But $D (-L)^{-1}=\int_0^\infty
DP_t \,dt$ according to \cite[Prop. 2.9.3]{IvanGioBook}.
Moreover,  
according\footnote{The statement of 
\cite[Prop. 5.1.5]{DavidEulaliaBook} is with $t^{-1/2}$ instead of $\frac{e^{-t}}{\sqrt{1-e^{-2t}}}$, but the given proof actually provides the estimate stated in (\ref{davideul}).}  to \cite[Prop. 5.1.5]{DavidEulaliaBook}, there exists $c_p>0$ such that, for any $F\in L^{p}(\Omega)$, 
\begin{equation}\label{davideul}
\|DP_tF\|_{L^{p}(\Omega,L^2(\R_+))}\leq c_p\,\frac{e^{-t}}{\sqrt{1-e^{-2t}}}\|F\|_{L^{p}(\Omega)}.
\end{equation}
 It follows from these two facts and the Minkowski inequality that the operator $D (-L)^{-1}$ is bounded from
 $L^{p}(\Omega)$ to $L^{p}(\Omega,L^2(\R_+))$. As a consequence, 
by iteration one has $$\sup_{\|h\|=1} \E\big[\| (D (-L)^{-1})^r(\varphi(W(h))) \|^{p}_{L^2(\R_+^r)}\big]
\leq c\,\E[|\varphi(N)|^{p}]<\infty$$
 for any $0\leq r\leq d$.

Let us finally consider the case $1\leq k\leq r$.
We can write, using among other the left-hand side of (\ref{meyer2}) and then its right-hand side,
\begin{eqnarray*}
&&\sup_{\|h\|=1} \E\big[\| D^{k}(D (-L)^{-1})^{r}(\varphi(W(h)))\|^{p}_{L^2(\R_+^{k+r})}\big]\\
&=& \sup_{\|h\|=1}\E\big[\| D^{k+1}(-L)^{-1}(D (-L)^{-1})^{r-1}(\varphi(W(h)))\|^{p}_{L^2(\R_+^{k+r})}\big]\\
&\leq&\frac{1}{c_1}\sup_{\|h\|=1}\E\big[\|(-L)^{\frac{k-1}2}(D (-L)^{-1})^{r-1}(\varphi(W(h)))\|^{p}_{L^2(\R_+^{r-1})}\big]\\
&\leq&\frac{c_{2}}{c_{1}}\sup_{\|h\|=1}\left(\E\big[\|D^{k-1}(D (-L)^{-1})^{r-1}\varphi(W(h))\|^{p}_{L^2(\R_+^{k+r-2})}\big]
\right.\\
&&\hskip4cm\left.
+\E\big[\|(D (-L)^{-1})^{r-1} \varphi(W(h))\|^{p}_{L^2(\R_+^{r-1})}\big]\right),
\end{eqnarray*}
which is finite by the induction property.
\qed

\section{Proof of Theorem~\ref{thm1}}
\label{sec3}

Since the point 1 (that is, convergence  of the finite-dimensional distributions) follows from the classical Breuer-Major theorem of \cite{BM}, 
let us only concentrate on the point 2.

We are thus left to show that the family  $(Y_n )_{n\ge 1} $ is tight in the Skorohod space ${\bf D}([0,1])$. 
Recall from \cite[Theorem 15.6]{Billingsley} that a sufficient condition for tightness in ${\bf D}([0,1])$  is the existence of $\gamma>0$ and $c>0$ such that, for all $n$,
\begin{equation}\label{titight}
\E[|Y_n(t)-Y_n(t_1)|^{1+\gamma}|Y_n(t_2)-Y_n(t)|^{1+\gamma}]\leq c\,(t_2-t_1)^{1+\gamma}, \quad 0\leq t_1\leq t\leq t_2\leq 1.
\end{equation}

We are not going to check (\ref{titight}) directly.
Instead, we shall use the following lemma, which is not stated in Billingsley book \cite{Billingsley} but has nevertheless become part of the folklore. For the sake of completeness, we give its proof.
\begin{lemma}\label{billylemma}
Fix $p>2$ and $c>0$.
If
\begin{equation}\label{tight}
\| Y_n(t) - Y_n(s) \|_{L^{p}(\Omega)}  \leq C  \left( \frac{ \lfloor nt \rfloor - \lfloor ns \rfloor } n \right)^{1/2},\quad s,t\in[0,1]
\end{equation}
for some $p>2$ and $C>0$
then (\ref{titight}) holds with $\gamma=\frac{p}2-1>0$ and $c=
3^{\frac{p}2}
C^{p}>0$.
\end{lemma}  
\noindent
{\it Proof}.
Suppose (\ref{tight}).
Using Cauchy-Schwarz, one has
\begin{eqnarray}
&&\E[|Y_n(t)-Y_n(t_1)|^{\frac{p}2}|Y_n(t_2)-Y_n(t)|^{\frac{p}2}]\notag\\
&\leq&
\| Y_n(t) - Y_n(t_1) \|_{L^{p}(\Omega)}^{\frac{p}2}
\| Y_n(t_2) - Y_n(t) \|_{L^{p}(\Omega)}^{\frac{p}2}\notag\\
&\leq&c^{p}
\left( \frac{ \lfloor nt \rfloor - \lfloor nt_1 \rfloor } n \right)^{\frac{p}4}
\left( \frac{ \lfloor nt_2 \rfloor - \lfloor nt \rfloor } n \right)^{\frac{p}4}.
\label{RHS}
\end{eqnarray}
If $\max(n(t-t_1),n(t_2-t))<\frac12$, then the quantity in (\ref{RHS}) is zero, and so (\ref{titight}) is verified.
If $n(t-t_1)\geq \frac12$, then 
$$
 \frac{ \lfloor nt \rfloor - \lfloor nt_1 \rfloor } n\leq  \frac{ nt - nt_1 +1 } n
\leq \frac{ nt - nt_1 +2n(t-t_1) } n \leq 3(t_2-t_1),
$$
whereas
$$
 \frac{ \lfloor nt_2 \rfloor - \lfloor nt \rfloor } n\leq  \frac{ nt_2 - nt +1 } n
\leq \frac{ nt_2 - nt +2n(t-t_1) } n \leq 3(t_2-t_1).
$$
Similar estimates hold if $n(t_2-t)\geq \frac12$.
So, if $\max(n(t-t_1),n(t_2-t))\geq \frac12$, then
the quantity in (\ref{RHS}) is bounded by $3^{\frac{p}2}c^{p}(t_2-t_1)^{\frac{p}2}$, and the proof of (\ref{titight}) is complete.
\qed  
  
\bigskip

We are now ready to proceed with the proof of point 2 in Theorem \ref{thm1}.
Combining the previous Lemma \ref{billylemma} with \cite[Theorem 15.6]{Billingsley}, we are left to show that
(\ref{tight}) is satisfied.

We can write
\begin{align*}
\| Y_n(t)- Y_n(s) \| _{L^{p}(\Omega)} &
= \frac1 {\sqrt{n}}  \left\| \sum_{i=\lfloor ns \rfloor} ^{\lfloor nt \rfloor-1} \varphi(X_i)  \right \| _{L^{p}(\Omega)} 
\!\!\!\!\!\!\!=  \frac1 {\sqrt{n}}   \left\| \sum_{i=\lfloor ns \rfloor} ^{\lfloor nt \rfloor-1}   \delta ^d \left( \varphi_d(X_i)    e_i^{\otimes d} \right)\right \| _{L^{p}(\Omega)} \mbox{by (\ref{varphir})}\\
&\le 
 c_p   \sum_{k=0}^d  \frac1 {\sqrt{n}}   \left\| \sum_{i=\lfloor ns \rfloor} ^{\lfloor nt \rfloor-1}   D^k \left( \varphi_d(X_i )   e_i ^{\otimes d} \right) \right \| _{L^{p}(\Omega; L^2(\R_+^{k+d}))}  \quad\mbox{by (\ref{meyer1})}\\
 & = c_p \sum_{k=0}^d  \left\|   \frac1 {n} \sum_{i,j=\lfloor ns \rfloor} ^{\lfloor nt \rfloor-1}   D^k (\varphi_d(X_i))  D^k  (\varphi_d(X_j))    \langle     e_i,  e_j \rangle_{L^2(\R_+)}  ^{d+k} \right \| _{L^{\frac{p}{2}}(\Omega; L^2(\R_+^k))}  ^{1/2}\\
 &=:  c_p   \sum_{k=0}^d  R_k.
\end{align*}
On the other hand,
\begin{eqnarray}
&&\sup_{i\in\Z} \| D^k (\varphi_d(X_i)) \|_{L^{p}(\Omega;L^2(\R_+^k))}  =\sup_{i\in\Z}  \E\big[
  \| D^k (\varphi_d(X_i)) \|^{p}_{L^2(\R_+^k))} 
 \big]^{\frac{1}{p}}\notag\\
 & =&\sup_{i\in\Z}  \E\big[
  \| D^k (
\langle (D (-L)^{-1})^d (\varphi(X_i)),e_i^{\otimes d}\rangle_{L^2(\R_+^d)}  
  )\|^{p}_{L^2(\R_+^k))}
 \big]^{\frac{1}{p}}\quad\mbox{by (\ref{varphid})}\notag\\
 &\leq &
\sup_{i\in\Z}  \E\big[
  \| D^k (D (-L)^{-1})^d (\varphi(X_i))\|^{p}_{L^2(\R_+^{k+d}))} 
 \big]^{\frac{1}{p}},\label{ck}
\end{eqnarray}
and (\ref{ck}) is finite thanks to Lemma \ref{lmcrucial}.

Recall from (\ref{ekrho}) that $\langle e_i,e_j\rangle_{L^2(\R_+)}=\rho(i-j)$.
Using Minkowski and H\"older inequalities, we can write, for any $0\leq k\leq d$,
\begin{align*}
R_k  & \le 
    \sup_{i\in\Z} \| D^k (\varphi_d(X_i)) \|_{L^{p}(\Omega;L^2(\R_+^k))}  \left(  \frac 1n  \sum_{i,j=\lfloor ns \rfloor} ^{\lfloor nt \rfloor-1}      | \rho(i-j)|^{d+k}    \right) ^{1/2}\\
    &\leq c_k\left(  \frac 1n  \sum_{i,j=\lfloor ns \rfloor} ^{\lfloor nt \rfloor-1}      | \rho(i-j)|^{d}    \right) ^{1/2}\quad\mbox{since $|\rho(k)|\leq 1$.}
\end{align*}

Finally,  the change of indices $(i,j) \to (i,j+h)$ leads to
 \[
 \frac 1n  \sum_{i,j=\lfloor ns \rfloor} ^{\lfloor nt \rfloor-1}     |    \rho(i-j)|^{d} 
 \le   C  \frac{ \lfloor nt \rfloor - \lfloor ns \rfloor } n \sum_{h\in \mathbb{Z}} | \rho(h)|^{d}  =  C\frac{ \lfloor nt \rfloor - \lfloor ns \rfloor } n,
 \]
 which provides the desired estimate (\ref{billylemma}) and concludes the proof of Theorem \ref{thm1}.
\qed

\section{Proof of Theorem~\ref{thm2}}
\label{sec4}

Since 
$
Z_n(t)=Y_n(t)+\frac{nt-\lfloor nt\rfloor}{\sqrt{n}}\,\varphi(X_{\lfloor nt\rfloor })
$
with
$\E\left[\left(\frac{nt-\lfloor nt\rfloor}{\sqrt{n}}\,\varphi(X_{\lfloor nt\rfloor })\right)^2\right]\leq \frac{1}{n}\,\|\varphi\|^2_{L^2(\gamma,\R)}\to 0
$,
point 1 (that is,
 convergence  of the finite-dimensional distribution of $Z_n$) follows again from the classical Breuer-Major theorem of  \cite{BM}.
 
Let us now turn to point 2.  
It remains to show that the family $(Z_n )_{n\ge 1} $ is tight in the space ${\bf C}([0,1])$. 
Recall from \cite[Theorem 12.3]{Billingsley} that
a sufficient condition for tightness in ${\bf C}([0,1])$ is this time the existence of $\gamma>0$ and $c>0$ such that, for all $n$,
\begin{equation}\label{billylemma2}
\| Z_n(t) - Z_n(s) \|_{L^{p}(\Omega)}  \leq c  |t-s|^{1/2},\quad s,t\in[0,1].
\end{equation}
Using the equivalent representation
$$
Z_n(t)=\frac{1}{\sqrt{n}}\int_0^{nt} \varphi(X_{\lfloor u\rfloor})du, 
$$
we can write
\begin{align*}
&\| Z_n(t)- Z_n(s) \| _{L^{p}(\Omega)} 
= \frac1 {\sqrt{n}}  \left\| \int_{ns}^{nt}\varphi(X_{\lfloor u\rfloor})du  \right \| _{L^{p}(\Omega)} \\
&=  \frac1 {\sqrt{n}}   \left\|\int_{ns}^{nt} \delta ^d \left( \varphi_d(X_{\lfloor u\rfloor}) e_{\lfloor u\rfloor}^{\otimes d} \right) du \right \| _{L^{p}(\Omega)} \quad\mbox{by (\ref{varphir})}\\
&\le 
 c_p   \sum_{k=0}^d  \frac1 {\sqrt{n}}   \left\| \int_{ns}^{nt} D^k \left( \varphi_d(X_{\lfloor u\rfloor} )   e_{\lfloor u\rfloor} ^{\otimes d} \right) du \right \| _{L^{p}(\Omega; L^2(\R_+^{k+d}))}  \mbox{by (\ref{meyer1})}\\
 & = c_p \sum_{k=0}^d  \left\|   \frac1 {n} \iint_{[ns,nt]^2}  D^k (\varphi_d(X_{\lfloor u\rfloor}))  D^k  (\varphi_d(X_{\lfloor v\rfloor}))    \langle     e_{\lfloor u\rfloor},  e_{\lfloor v\rfloor} \rangle_{L^2(\R_+)}  ^{d+k}  dudv\right \| _{L^{\frac{p}{2}}(\Omega; L^2(\R_+^k))}  ^{1/2}\\
 &=:  c_p   \sum_{k=0}^d  R_k.
\end{align*}

Minkowski and H\"older inequalities yield, for any $0\leq k\leq d$,
\begin{align*}
R_k  & \le 
    \sup_{u\in\R_+} \| D^k (\varphi_d(X_{\lfloor u\rfloor})) \|_{L^{p}(\Omega;L^2(\R_+^k))}  \left(  \frac 1n  \iint_{[ns,nt]^2} | \rho({\lfloor u\rfloor}-{\lfloor v\rfloor})|^{d+k}  dudv  \right) ^{1/2}\\
    &= c_k\left(  \frac 1n  \iint_{[ns,nt]^2} | \rho({\lfloor u\rfloor}-{\lfloor v\rfloor})|^{d+k}  dudv  \right) ^{1/2},
\end{align*}
with $c_k$ finite by (\ref{ck}).

Finally,  since $|\rho(k)|\leq 1$ for all $k$,
 \begin{eqnarray*}
 \frac 1n  \iint_{[ns,nt]^2} | \rho({\lfloor u\rfloor}-{\lfloor v\rfloor})|^{d+k}  dudv  
&\leq&\frac1n\int_{ns}^{nt} \left( \int_{ns-v}^{nt-v}| \rho(\lfloor x+v\rfloor-\lfloor v\rfloor)|^{d} dx\right)dv\\ 
  &\le& \frac1n\int_{ns}^{nt} \sum_{j\in\Z}|\rho(j)|^{d}dv
  =  |t-s|\, \sum_{j\in\Z}|\rho(j)|^{d},
 \end{eqnarray*}
 which provides the desired estimate (\ref{billylemma2}) and concludes the proof of Theorem \ref{thm2}.
\qed

\section{An extension of Theorem \ref{thm1}}

In this section, our aim is to show that the method we have employed
for the proofs of Theorems \ref{thm1} and \ref{thm2} can be easily extended to deal with the case
where $\sum_{|j|\leq n}|\rho(j)|^d$ diverges as a slowly varying function when $n\to\infty$.
Instead of stating such a result at a great level of generality, to avoid too much technicalities we prefer to illustrate what happens in a guiding example and
only in the setting of Theorem \ref{thm1}. The same extension for Theorem \ref{thm2} would follow similar lines; details are left to the interested reader as an exercise.

Consider the fractional Gaussian noise 
$X_k=B_{k+1}-B_k$ associated with
a fractional Brownian motion $B$ of Hurst index $H\in (0,1)$; in this case,
$\rho(k)=\frac12\big(|k+1|^{2H}+|k-1|^{2H}-2|k|^{2H}\big)$.
Also, consider a function  $\varphi \in L^2(\R,\gamma)$ with expansion (\ref{hermite}) and Hermite rank $d\ge 1$.
Finally, recall $Y_n$ from (\ref{yn}), let $W=\{W_t\}_{t\in [0,1]}$ be a Brownian motion and let $\sigma ^2$ be defined in (\ref{bm.sig}). 

Since $\rho(k)\sim c|k|^{2H-2}$ (where
$c$ is an explicit constant whose value is useless), in the case where $H\in(0,1-\frac{1}{2d})$ one can apply Breuer-Major theorem of \cite{BM} to deduce that
$Y_n\overset{\rm f.d.d}{\to} \sigma W$. If moreover $\varphi \in L^p(\R,\gamma)$ for some $p>2$, then $Y_n\overset{{\bf D}([0,1])}{\to} \sigma W$ thanks to our Theorem \ref{thm1}.

In contrast, when  $H\in(1-\frac{1}{2d},1)$ Taqqu \cite[Theorem 5.6]{Taqqu}  has shown in the seventies that 
$$n^{d(1-H)-\frac12}Y_n\overset{{\bf D}([0,1])}{\to} Y_\infty,$$ 
where $Y_\infty$ stands for the Hermite process of index $d$. Here, note that no additional integrability condition on $\varphi$ is required for the convergence to hold in ${\bf D}([0,1])$; indeed, since the limiting process $Y_\infty$ is $\alpha$-H\"older continuous with $\alpha$ \emph{strictly} greater than $\frac12$, it is enough to bound $\|Y_n(t)-Y_n(s)\|_{L^2(\Omega)}$ (and not $\|Y_n(t)-Y_n(s)\|_{L^p(\Omega)}$ with $p>2$) to get the tightness, so classical and easy calculations are enough to conclude.

What about the critical case $H=1-\frac{1}{2d}$?
In this case, $\rho(k)\sim c|k|^{-\frac1d}$ and so $\sum_{k\in\Z}|\rho(k)|=+\infty$. Nevertheless, since the divergence of the series is slow, the fluctuations are still Brownian after proper normalisation.
More precisely, it is shown in \cite{BM} that
$\frac{Y_n}{\sqrt{\log n}}\overset{\rm f.d.d}{\to} \sigma W$, 
with $\sigma^2=2d!\big(\frac{(2d-1)(d-1)}{2d^2}\big)^d$.
As far as the convergence in ${\bf D}([0,1])$ is concerned, a slight extension of our method leads to the following result.

\begin{theorem} \label{thm3}
Consider a function  $\varphi \in L^2(\R,\gamma)$ with expansion (\ref{hermite}) and Hermite rank $d\ge 1$.
Let $X=\{ X_n\}_{n\in \mathbb{Z}} $ be the fractional Gaussian noise
of index $H=1-\frac1{2d}$, that is, $X$ is a mean-zero Gaussian stationary sequence
with convariance function
$$
\E[X_nX_{n+k}]=\rho(k)=\frac12\big(|k+1|^{2-\frac1d}+|k-1|^{2-\frac1d}-2|k|^{2-\frac1d}\big).
$$ 
Finally, recall $Y_n$ from (\ref{yn}) , let $W=\{W_t\}_{t\in [0,1]}$ be a Brownian motion and let $\sigma$ be given by
$
\sigma^2=2d!\big(\frac{(2d-1)(d-1)}{2d^2}\big)^d.
$
Then, as $n\to\infty$,
\begin{enumerate}
\item The finite-dimensional distributions of $\frac{Y_n}{\sqrt{\log n}}$ converge to those of $\sigma W$;
\item If  $\varphi \in L^{p}(\R,\gamma)$  for some $p>2$, then
$\frac{Y_n}{\sqrt{\log n}}$ converges in   law to $\sigma W$ in ${\bf D}([0,1])$ 
endowed with the Skorohod topology.
\end{enumerate}
\end{theorem}
\noindent
{\it Proof}. Point 1 follows from Breuer and Major \cite{BM}.
  Combining Lemma \ref{billylemma} with \cite[Theorem 15.6]{Billingsley} (for  $\frac{Y_n}{\sqrt{\log n}}$ instead of $Y_n$), to prove point 2 it is enough to show that
(\ref{tight}) holds true.

We can write
\begin{align*}
&\left\| \frac{Y_n(t)}{\sqrt{\log n}} - \frac{Y_n(s)}{\sqrt{\log n}} \right\| _{L^{p}(\Omega)} \\
&
=  \frac1 {\sqrt{n\log n}}   \left\| \sum_{i=\lfloor ns \rfloor} ^{\lfloor nt \rfloor-1}   \delta ^d \left( \varphi_d(X_i)    e_i^{\otimes d} \right)\right \| _{L^{p}(\Omega)} \mbox{by (\ref{varphir})}\\
&\le 
 c_p   \sum_{k=0}^d  \frac1 {\sqrt{n\log n}}   \left\| \sum_{i=\lfloor ns \rfloor} ^{\lfloor nt \rfloor-1}   D^k \left( \varphi_d(X_i )   e_i ^{\otimes d} \right) \right \| _{L^{p}(\Omega; L^2(\R_+^{k+d}))}  \quad\mbox{by (\ref{meyer1})}\\
 & = c_p \sum_{k=0}^d  \left\|   \frac1 {n\log n} \sum_{i,j=\lfloor ns \rfloor} ^{\lfloor nt \rfloor-1}   D^k (\varphi_d(X_i))  D^k  (\varphi_d(X_j))    \langle     e_i,  e_j \rangle_{L^2(\R_+)}  ^{d+k} \right \| _{L^{\frac{p}{2}}(\Omega; L^2(\R_+^k))}  ^{1/2}\\
 &\leq   c_p   \sum_{k=0}^d  
 \sup_{i\in\Z} \| D^k (\varphi_d(X_i)) \|_{L^{p}(\Omega;L^2(\R_+^k))}  \left(  \frac 1{n\log n}  \sum_{i,j=\lfloor ns \rfloor} ^{\lfloor nt \rfloor-1}      | \rho(i-j)|^{d+k}    \right) ^{1/2}\\
    &\leq  c  \left(  \frac 1{n\log n}  \sum_{i,j=\lfloor ns \rfloor} ^{\lfloor nt \rfloor-1}      | \rho(i-j)|^{d}    \right) ^{1/2}\quad\mbox{since $|\rho(k)|\leq 1$ and using (\ref{ck}).}
\end{align*}

Finally,  the change of indices $(i,j) \to (i,j+h)$
and the fact that $|\rho(r)|^d\sim c|r|^{-1}$ as $|r|\to\infty$ leads to
 \[
 \frac 1{n\log n}  \sum_{i,j=\lfloor ns \rfloor} ^{\lfloor nt \rfloor-1}     |    \rho(i-j)|^{d} 
 \le   C\frac{ \lfloor nt \rfloor - \lfloor ns \rfloor } n,
 \]
 which provides the desired estimate (\ref{billylemma}) and concludes the proof of Theorem \ref{thm1}.
\qed

\end{document}